\newif\ifpreprint
\preprinttrue   % <-- Change to \preprintfalse for anonymous submission

\ifpreprint
  \documentclass[pdflatex,sn-mathphys-num]{sn-jnl}
\else
  \documentclass[referee,pdflatex,sn-mathphys-num]{sn-jnl}
\fi

\usepackage{graphicx}
\usepackage{amsmath,amssymb,amsfonts}
\usepackage{booktabs}
\usepackage{setspace}

\ifpreprint
  \renewcommand{\normalsize}{\fontsize{11}{13}\selectfont}  % 11pt font
  \normalsize
\fi

\begin{document}

%%======================================================================
%% TITLE
%%======================================================================
\title[A dynamical model of the U.S. mathematics graduate degree pipeline]{A Dynamical Model of the U.S. Mathematics Graduate Degree Pipeline}

%%======================================================================
%% AUTHORS AND AFFILIATIONS
%%======================================================================
\ifpreprint
  % DE-ANONYMIZED VERSION (for preprint)
  \author*[1,2,3]{\fnm{Chad M.} \sur{Topaz}}\email{chad@qsideinstitute.org}
  \author[4]{\fnm{Oluwatosin} \sur{Babasola}}
  \author[5]{\fnm{Ron} \sur{Buckmire}}
  \author[6]{\fnm{Daozhou} \sur{Gao}}
  \author[7]{\fnm{Maila} \sur{Hallare}}
  \author[8]{\fnm{Olaniyi} \sur{Iyiola}}
  \author[9]{\fnm{Deanna} \sur{Needell}}
  \author[10]{\fnm{Andr\'es R.} \sur{Vindas-Mel\'endez}}

  \affil[1]{\orgname{Williams College}, \city{Williamstown}, \state{MA}, \country{USA}}
  \affil[2]{\orgname{QSIDE Institute}, \city{Williamstown}, \state{MA}, \country{USA}}
  \affil[3]{\orgname{University of Colorado}, \city{Boulder}, \state{CO}, \country{USA}}
  \affil[4]{\orgname{University of Georgia}, \city{Athens}, \state{GA}, \country{USA}}
  \affil[5]{\orgname{Marist University}, \city{Poughkeepsie}, \state{NY}, \country{USA}}
  \affil[6]{\orgname{Cleveland State University}, \city{Cleveland}, \state{OH}, \country{USA}}
  \affil[7]{\orgname{United States Air Force Academy}, \state{CO}, \country{USA}}
  \affil[8]{\orgname{Morgan State University}, \city{Baltimore}, \state{MD}, \country{USA}}
  \affil[9]{\orgname{University of California}, \city{Los Angeles}, \state{CA}, \country{USA}}
  \affil[10]{\orgname{Harvey Mudd College}, \city{Claremont}, \state{CA}, \country{USA}}
\else
  % ANONYMOUS VERSION (for journal submission)
  \author{Anonymous Author(s)}
  \affil{Anonymous Affiliation(s)}
\fi

%%======================================================================
%% ABSTRACT (150-250 words per La Matematica guidelines)
%%======================================================================
\abstract{We present a latent-stock compartmental framework for modeling degree production systems when only completion flows, rather than enrollments, are observed. Applied to U.S.\ mathematics degrees from 1969 to 2017, the model treats master's and PhD populations as latent compartments---unobserved state variables that are inferred indirectly because they generate the observed completion flows---with time-varying routing fractions and completion hazards. Using information-criterion model comparison across a grid of specifications, we find strong support for smooth nonlinear time variation in routing fractions and hazards, while models with explicit international forcing are disfavored. The preferred model achieves a log-scale root mean squared error of approximately 0.036, corresponding to a typical multiplicative error of about 4\% in fitted degree counts, and highlights key structural shifts in the graduate pipeline: the master's pathway became increasingly central to PhD production through the late twentieth century before weakening, while direct bachelor's-to-PhD entry remained small but persistent. Estimated completion hazards for both degrees rise over time, indicating faster effective turnover in the graduate compartments. Methodologically, our main contribution is a latent stock dynamical approach that recasts linked degreecompletion time series as a coherent stock-flow system when intermediate enrollments are unobserved, making explicit both what features of pipeline dynamics are identifiable from completion data alone and what limitations such data impose.}

%%======================================================================
%% KEYWORDS (4-6 per La Matematica guidelines)
%%======================================================================

%% MSC Classification codes (add to separate title page or as journal requires):
%% 97B40: Higher education
%% 91D30: Social networks; opinion dynamics
%% 37N40: Dynamical systems in optimization and economics
%% 62P25: Applications of statistics to social sciences

\maketitle

\section{Introduction}

This paper develops a latent-stock compartmental framework for modeling degree-production systems when only aggregate completion flows, rather than enrollments, are observed. We apply this framework to the U.S.\ mathematics training pipeline, connecting the annual production of bachelor's, master's, and PhD degrees over a span of nearly five decades. 
Our aim is to understand how changes in bachelor's degree production propagate forward through the system to influence master's and PhD completions, using only publicly available degree counts in a reduced-form descriptive framework.

Understanding these dynamics matters for both substantive and methodological reasons. 
Substantively, graduate education operates as a pipeline with long delays: shifts in degree production today can shape advanced training, labor supply, and the academic workforce many years later. 
Structural changes in how students move through the pipeline may therefore have long-lasting consequences that are not apparent from contemporaneous degree counts alone.
Methodologically, however, the quantities that would most naturally describe such a pipeline in this context---annual enrollments, i.e.\ intermediate stocks---are typically unobserved, leaving only aggregate completion flows available for analysis. 
This combination makes it difficult to distinguish short-run fluctuations from long-run structural change without an explicit dynamical framework.

\begin{figure}[t!]
\centering
\includegraphics[width=\textwidth]{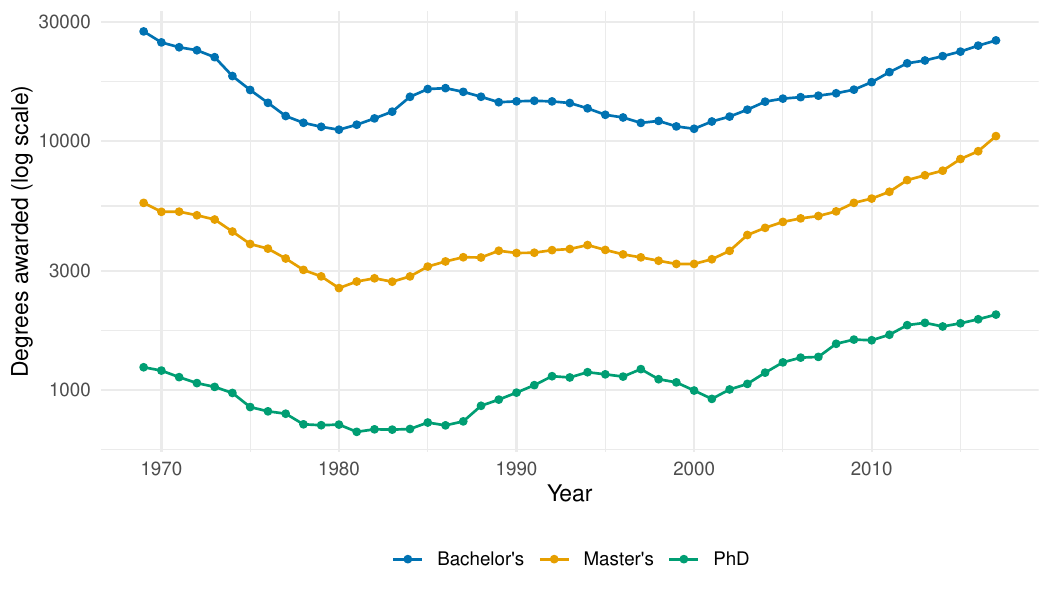}
\caption{Log-scale annual observations of bachelor's, master's, and PhD degrees in mathematics awarded in the United States for 49 observation years between 1969 and 2017. Earlier nonannual observations prior to 1969 are omitted from the figure and from the model estimation. 
Source: National Center for Education Statistics (NCES), U.S.\ Department of Education.}
\label{fig:raw}
\end{figure}

The key empirical constraint shaping our analysis of this problem is that national data sources such as the National Center for Education Statistics at the United States Department of Education report only the \emph{flows} of degrees awarded per year, not the underlying enrollments. 
For the modeling work below we use 49 annual observations between 1969 and 2017, during which we observe three series: the number of bachelor's degrees $b(t)$, master's degrees $m(t)$, and PhD degrees $p(t)$ awarded in mathematics. 
These are completion counts---snapshots of how many students finished each level of training in a given year. 
We do not have access to data for how many students are currently enrolled at each level, how long they have been enrolled, or which specific individuals transition from one level to the next. 
Figure~\ref{fig:raw} plots these three series on a log scale. 
All exhibit a decline in the 1970s, recovery through the 1980s and 1990s, and rapid growth in the 2000s and 2010s. 
Their co-movement over time motivates treating them as components of a linked dynamical system rather than as independent sequences.

This data structure poses a fundamental modeling challenge. 
If enrollments and retention or attrition dynamics were observed, one might attempt to track cohorts of students as they progress through the pipeline. 
In the absence of enrollment data, however, analysis must instead proceed from aggregate flows alone.
The central question then becomes: \emph{what is the simplest dynamical system that takes bachelor's completions as an input and produces master's and PhD completions as outputs in a way that is consistent with the historical record?}

A natural framework for addressing this question is a compartmental model, in which populations move through a sequence of stages governed by stock--flow relationships---a modeling approach with a long history in epidemiology, demography, and system dynamics (see, e.g., \cite{Hethcote2000,Forrester1961,Sterman2000}). 
The appeal of such models lies in their ability to connect stocks (how many individuals occupy each stage) to flows (how many individuals move between stages per unit time). 
In our setting, the relevant structure is that of a linear compartmental system with immigration: each compartment receives new entrants and experiences exits at rates that may vary over time, and the PhD compartment receives an additional inflow from an upstream compartment. 
Because the available data are annual, we work directly with a discrete-time formulation, aligning the model with the temporal resolution of the observations and avoiding additional numerical discretization.

Applying compartmental models in this context raises two challenges. First, such models are naturally expressed in terms of compartment populations, yet the populations of interest here are entirely unobserved. We observe only exits---the degrees awarded---not the stocks of students currently enrolled. The intermediate populations must therefore be treated as latent variables, inferred indirectly from the degree flows they generate. Second, standard compartmental models often assume constant transition rates, whereas the structure of graduate education has evolved substantially over the five decades we study. As we will show, a model with constant parameters cannot adequately capture the observed degree series. Our strategy is therefore to introduce \emph{latent compartments} representing the unobserved populations of students at the master's and PhD levels. These compartments accumulate students over time as new entrants arrive from earlier stages of the pipeline and release students as degrees are completed. The observable degree flows then emerge as outputs of this hidden dynamical process. By fitting the model parameters to match the observed degree series, we can infer how effective routing behavior and aggregate turnover have evolved over the past five decades, even though the underlying populations are never directly observed. The resulting model is a reduced-form dynamical approximation: it is designed to summarize aggregate patterns in the data rather than to provide a mechanistic or causal account of individual student decisions.

While related challenges arise in epidemiology, demography, labor economics, and system dynamics, existing applications in those fields typically rely on auxiliary stock information, partial enrollment data, or equilibrium assumptions, rather than inference from exit flows alone. By contrast, our setting requires inference from completion flows alone. In this sense, our estimation problem can be viewed as an \emph{inverse stock--flow problem}, in which a parsimonious latent dynamical system is inferred from observed exits. The framework developed here therefore contributes a general approach to studying stock--flow systems observed only through their exits, while also making explicit which features of long-run dynamics are identifiable from such data and which are intrinsically weakly identified.

The remainder of this paper proceeds as follows. We first develop a two-compartment model with constant parameters, establishing the basic logic of how compartments, inflows, and outflows relate to one another. We then generalize the model to allow parameters to vary smoothly over time, which turns out to be essential for capturing the long-run evolution of the degree series. After describing how we estimate the parameters from the observed data, we discuss model selection and present the fitted results, including uncertainty quantification based on the curvature of the loss function. Finally, we consider whether explicitly modeling international student inflows improves the model and discuss how the findings should be interpreted.

\section{Background}

A growing body of work on the academic workforce treats higher education not as a collection of static snapshots but as a dynamic population in which individuals move through structured stages with temporal lags, demographic constraints, and bottlenecks. In mathematics, long-run patterns of training and hiring evolve cumulatively over time: analyses of faculty hiring networks show that prestige, gender, and subfield jointly shape year-to-year hiring flows in ways that cannot be understood cross-sectionally \cite{FitzGeraldEtAl2023}. Related demographic studies demonstrate that even large increases in women's representation among new PhD cohorts translate only slowly into changes in faculty composition, because age distributions, turnover rates, and limited replacement opportunities constrain the speed of adjustment \cite{HargensLong2002}. Such findings reflect a broader theoretical view that organizational populations exhibit strong inertia, with incumbency and slow exit dynamics limiting responsiveness to external shocks \cite{HannanFreeman1977}.

Building on these foundations, several lines of work develop explicit dynamical models of academic populations. Stage- or age-structured demographic models show that faculty representation changes only gradually even under substantially altered hiring or retention policies \cite{MarschkeEtAl2007}. Stochastic update models calibrated to census-level data confirm that hiring and attrition make distinct contributions to shifts in faculty gender composition, and that counterfactual interventions---such as equalizing retention risks or modifying hiring proportions---produce only modest short-term effects because of demographic inertia \cite{LaBergeEtAl2024}. Pipeline-oriented studies extend this logic beyond the faculty stage: comparisons of stage-neutral and stage-specific transition systems isolate where disproportionate losses occur for different demographic groups \cite{ShawStanton2012,ShawEtAl2021}. Related work on the postdoctoral workforce demonstrates how seemingly small differences in transition rates accumulate as cohorts move through the academic pipeline, reinforcing or amplifying disparities \cite{GhaffarzadeganHawleyDesai2014}. Collectively, these results highlight the importance of modeling transition flows explicitly rather than inferring long-run outcomes from endpoint comparisons alone.

While much of this literature focuses on downstream career stages, degree production itself constitutes the upstream supply process that constrains all subsequent academic labor markets. Shifts in the volume, routing, and timing of degree completions shape postdoctoral and faculty dynamics long before they appear in hiring or retention data. Understanding degree-production dynamics is therefore essential for interpreting long-run changes in the academic workforce, rather than merely a preliminary step toward modeling later career stages.

A complementary set of approaches frames academic population dynamics through system-level constructs such as reproduction numbers, queueing relationships, and branching processes. Analyses of PhD production rates show that typical faculty members ``reproduce'' far more doctorates over a career than can be absorbed into the tenure-track system, implying structural oversupply unless faculty positions expand exponentially \cite{LarsonGhaffarzadeganXue2014}. Extensions of this framework argue that oversupply at the PhD-to-faculty transition can coexist with shortages elsewhere in the STEM labor market \cite{XueLarson2015}. Queueing-based and stock-and-flow models make explicit how delays and turnover shape long-run system behavior, demonstrating that lengthening academic careers mechanically reduces the number of available entry-level faculty positions over extended horizons \cite{LarsonGomezDiaz2012}. System-dynamics models at the institutional level further illustrate how enrollment, staffing, and resource constraints interact through low-dimensional flow structures \cite{Zaini2017}, while population-process models show that branching-like reproduction and lifetime dynamics can capture long-run behaviors in scholarly ecosystems \cite{WuVenkatramananChiu2015}.

Across these literatures, a consistent insight emerges: academic systems are inherently dynamical, and their long-run evolution depends on explicitly modeling flows, transition rates, and structural constraints. Existing research, however, focuses predominantly on faculty hiring, retention, postdoctoral transitions, or career-stage dynamics; comparatively little work models the degree-production process itself as a linked dynamical system using degree counts alone. Moreover, while related challenges arise in demography, epidemiology, labor economics, and system dynamics, most applications in those fields rely on observed or partially observed stocks, auxiliary enrollment data, or equilibrium assumptions. By contrast, the setting considered here requires inference when intermediate enrollments are entirely unobserved and only aggregate completion flows are available.

To our knowledge, there is limited prior work that links annual bachelor's, master's, and PhD degree counts through latent graduate compartments inferred solely from completions, while allowing routing behavior and aggregate turnover to vary smoothly over multi-decade horizons. The mathematics training system---characterized by stable degree definitions, long observational records, and clearly delineated transitions---is particularly well suited for such an analysis. The compartmental framework developed below recasts three separate degree-count time series as a single, dynamically coherent system and clarifies the extent to which long-run changes in pipeline routing and turnover can be inferred from degree counts alone.

\section{Modeling}

This section develops a two-compartment dynamical model that links bachelor's, master's, and PhD degree production in mathematics. We begin with a baseline specification in which all parameters are constant over time, and then extend the model to allow smooth time variation in the branching fractions and completion hazards so that it can accommodate long-run changes in the graduate training pipeline. For reference, Table~\ref{tab:notation} summarizes all notation used in the modeling framework and throughout the remainder of the paper.

\begin{table}[t!]
\centering
\renewcommand{\arraystretch}{1.15}
\caption{Notation for the degree-production model. Time $t$ indexes academic years (annual resolution). Degree counts are annual completions (flows). Graduate compartments are latent (unobserved) stocks inferred from completion flows.}
\label{tab:notation}
\begin{tabular}{
  p{0.28\textwidth}
  @{\hspace{0.03\textwidth}}
  p{0.63\textwidth}
}\toprule
\textbf{Symbol} & \textbf{Meaning} \\
\midrule
$t$ & Academic year index. \\
$t_{\min},\, t_{\max}$ & Earliest and latest observation years (1969, 2017). \\
$t_{\mathrm{mid}}$ & Midpoint year $(t_{\min}+t_{\max})/2$ (1993). \\
$s(t)$ & Rescaled time index mapping $t$ to approximately $[-1,1]$. \\
\midrule
$b(t)$ & Observed bachelor's degree completions in mathematics. \\
$m(t)$ & Observed master's degree completions in mathematics. \\
$p(t)$ & Observed PhD degree completions in mathematics. \\
\midrule
$M(t)$ & Latent master's compartment stock (effective degree-producing population). \\
$P(t)$ & Latent PhD compartment stock (effective degree-producing population). \\
\midrule
$\hat{m}(t)$ & Model-implied master's completions; $\hat{m}(t)=\gamma_M(t)M(t)$. \\
$\hat{p}(t)$ & Model-implied PhD completions; $\hat{p}(t)=\gamma_P(t)P(t)$. \\
\midrule
$\rho_{BM},\,\rho_{BP},\,\rho_{MP}$ & Effective branching fractions (constant-parameter model). \\
$\rho_{BM}(t),\,\rho_{BP}(t),\,\rho_{MP}(t)$ & Time-varying effective branching fractions. \\
\midrule
$\gamma_M,\,\gamma_P$ & Completion hazards (annual completion fractions; constant model). \\
$\gamma_M(t),\,\gamma_P(t)$ & Time-varying completion hazards. \\
$1/\gamma(t)$ & Effective duration (years) under memoryless turnover approximation. \\
\midrule
$\mathrm{logit}(x)$ & Logit transform $\log\!\big(x/(1-x)\big)$, $x\in(0,1)$. \\
$\mathrm{logit}^{-1}(y)$ & Logistic inverse $(1+e^{-y})^{-1}$, $y\in\mathbb{R}$. \\
\bottomrule
\end{tabular}
\end{table}

\subsection{Model with Constant Parameters}

We begin by developing the simplest version of our model, in which all parameters are constant over time. This allows us to establish the basic structure and introduce key concepts before adding the complications needed to fit the observed data.

We index time in one-year increments and write $t$ for the academic year. Let $M(t)$ denote the latent population in the master's compartment at time $t$, and let $P(t)$ denote the latent population in the PhD compartment. We do not observe $M(t)$ or $P(t)$ directly; they are \emph{latent} state variables inferred from observed degree completions. The quantities we do observe are the annual degree completions: $b(t)$ bachelor's degrees, $m(t)$ master's degrees, and $p(t)$ PhD degrees.

The central modeling assumption is that each compartment behaves like a reservoir with inflows and outflows over each one-year interval. Students enter the master's compartment when they transition into the effective master's stock that generates degree completions under the model, and they leave when they complete their master's degree. Similarly, students enter the PhD compartment when they transition into the effective PhD stock and leave upon completion. This language reflects the fact that the model represents effective degree-producing populations rather than literal enrollments.

We do not include an explicit attrition (non-completion) outflow from the graduate compartments; students who leave graduate study without completing are not modeled as a distinct flow. Accordingly, the latent stocks $M(t)$ and $P(t)$ should be interpreted as \emph{effective degree-producing stocks} rather than literal enrollments. In this reduced-form representation, non-completion exits are not separately identified; their net effect is absorbed into the fitted effective routing fractions and effective turnover implied by the hazards, rather than appearing as an explicit outflow.

We treat the bachelor's completion series $b(t)$ as an exogenous input to the model. This is a deliberate simplification: we do not attempt to model the dynamics of undergraduate enrollment or completion, but instead take the observed bachelor's degree output as given and ask how it feeds into the downstream stages. This choice reflects both the data we have available and a desire to focus on the graduate training pipeline.

When a student completes a bachelor's degree in mathematics, several outcomes are possible: they might enter a master's program in mathematics, enter a PhD program directly, or leave the mathematics pipeline entirely. We summarize the aggregate pattern of these transitions using \emph{branching fractions}.

Let $\rho_{BM}$ denote the effective fraction of bachelor's completions that feed into the master's compartment, and let $\rho_{BP}$ denote the effective fraction that feed directly into the PhD compartment. The word \emph{effective} is important here and deserves emphasis. Throughout this paper, we use \emph{effective branching fraction} and \emph{effective hazard} to mean aggregate parameters calibrated to match observed flows, as distinct from individual-level transition probabilities that might be measured from student records. These effective parameters incorporate heterogeneity, selection, and timing effects that we do not model explicitly because doing so would require individual-level or enrollment data that are not available in the national degree series.

These branching fractions are not probabilities measured at the individual level, but rather aggregate parameters that describe how the bachelor's completion flow contributes to downstream transitions in our reduced-form system. They incorporate the net effects of domestic and international students, of students who enter immediately versus those who delay, and of various institutional pathways that we do not model explicitly. We do not model delayed entry or leave-and-return behavior with explicit lag structure; their net effect is absorbed into the effective routing fractions and hazards. Importantly, we estimate $\rho_{BM}$ and $\rho_{BP}$ independently and do not impose a simplex constraint. Imposing a simplex constraint (e.g., $\rho_{BM} + \rho_{BP} = 1$) would force the model to allocate essentially all bachelor's completions to modeled graduate-entry pathways. Because many bachelor's recipients do not enter U.S.\ mathematics graduate programs, or do so only after substantial delays, we allow $\rho_{BM} + \rho_{BP} < 1$ and interpret the remainder as an implicit exit fraction from the mathematics pipeline. In our fitted time-varying model below, we find empirically that $\rho_{BM}(t) + \rho_{BP}(t) < 1$ for all $t$, reflecting the inferred dynamics rather than an imposed algebraic restriction.

Similarly, let $\rho_{MP}$ denote the effective fraction of master's completions that feed into the PhD compartment. It is important to note that $\rho_{MP}$ describes the effective fraction of \emph{completed} master's degrees that eventually contribute to PhD enrollment, not the share of master's enrollees who go on to doctoral study. The distinction matters because the population of master's completers differs from the population of master's enrollees, and because some students who eventually pursue PhDs do so years after completing their master's degree. We do not model such delays explicitly; instead, the aggregate impact of these lags is absorbed into the reduced-form fraction $\rho_{MP}$.

With these definitions, the inflow to the master's compartment during year $t$ is $\rho_{BM} \, b(t)$. The inflow to the PhD compartment has two sources: direct entry from bachelor's completions, contributing $\rho_{BP} \, b(t)$, and entry from master's completions, contributing $\rho_{MP} \, \hat{m}(t)$, where $\hat{m}(t)$ denotes model-implied master's completions.

To model the outflows from each compartment, we introduce annual \emph{completion hazards}, also called per-capita completion fractions. Let $\gamma_M$ denote the master's completion hazard and $\gamma_P$ the PhD completion hazard. The model-implied number of master's degrees completed in year $t$ is
\begin{equation}
\hat{m}(t) = \gamma_M \, M(t),
\end{equation}
and the model-implied number of PhD degrees completed in year $t$ is
\begin{equation}
\hat{p}(t) = \gamma_P \, P(t).
\end{equation}
For constant hazards, the reciprocals $1/\gamma_M$ and $1/\gamma_P$ can be interpreted as steady-state average durations spent in each compartment under a memoryless exit approximation. When hazards vary over time, it is more appropriate to interpret $1/\gamma(t)$ as an instantaneous effective duration rather than a literal cohort-average time-to-degree.

As with the branching fractions, we emphasize that $\gamma_M$ and $\gamma_P$ are effective parameters for the aggregate system. They summarize how quickly the effective degree-producing stocks turn over in terms of degree production, but they should not be interpreted as direct measurements of time-to-degree for any particular subgroup of students.

The evolution of the two compartments is described by discrete-time update equations at annual resolution:
\begin{equation}
M(t+1) = M(t) + \rho_{BM} \, b(t) - \gamma_M \, M(t),
\end{equation}
\begin{equation}
P(t+1) = P(t) + \rho_{BP} \, b(t) + \rho_{MP} \, \hat{m}(t) - \gamma_P \, P(t).
\end{equation}
Here the two compartments are dynamically coupled because the PhD inflow term $\rho_{MP}\,\hat m(t)$ depends on the master's stock via $\hat m(t)=\gamma_M(t)M(t)$. Throughout estimation we constrain parameters so that the system remains in the nonnegative orthant: completion hazards satisfy $\gamma_M(t), \gamma_P(t)\in(0,1)$, effective branching fractions satisfy $\rho_{BM}(t), \rho_{BP}(t), \rho_{MP}(t)\ge 0$, the observed input $b(t)$ is nonnegative, and initial conditions satisfy $M(t_{\min})>0$ and $P(t_{\min})>0$. Under these conditions the discrete-time updates preserve nonnegativity of the latent stocks for all $t$. We interpret our equations so that completions observed in year $t$ contribute to enrollment at the start of year $t+1$, consistent with the academic calendar in which many bachelor's graduates enter graduate programs in the subsequent academic year.

These are linear first-order difference equations with time-varying forcing through the exogenous series $b(t)$ and the derived series $\hat{m}(t)$. Linearity in the state variables $M$ and $P$ ensures a unique state trajectory for any given initial conditions and parameter values, simplifying estimation and interpretation.

Given the observed bachelor's series $b(t)$, initial conditions $M(t_{\min})$ and $P(t_{\min})$, and candidate parameter values, these equations determine the trajectories of $M(t)$ and $P(t)$, and hence the model-implied completions
\begin{equation}
\hat{m}(t) = \gamma_M\, M(t), \qquad \hat{p}(t) = \gamma_P\, P(t).
\end{equation}
In this discrete-time setting, $\gamma_M$ and $\gamma_P$ are per-year completion fractions rather than continuous-time hazard rates. When $\gamma$ is approximately stable over a period, $1/\gamma$ provides an interpretable effective duration in years under the reduced-form turnover approximation. Together, the parameters $\rho_{BM}$, $\rho_{BP}$, $\rho_{MP}$, $\gamma_M$, and $\gamma_P$ govern how bachelor's output feeds forward through the pipeline and how rapidly the graduate compartments turn over.

Because the stocks are latent, the model requires a convention for initializing $M(t_{\min})$ and $P(t_{\min})$. In our empirical implementation, we choose initial conditions that match the first observed master's and PhD completion counts exactly given the fitted hazards; the specific convention is described in Section~\ref{sec:estimation}.

\subsection{Time-Varying Parameters}

The constant-parameter model captures the basic logic of the training pipeline, but it cannot fit the observed data well. Over the decades spanned by our sample, the structure of mathematics graduate education changed substantially, as reflected directly in the long-run evolution of the degree series themselves. To accommodate these long-run changes, we make a key modeling choice: we allow the branching fractions and completion hazards to vary smoothly over time.

Introducing time-varying parameters creates both opportunities and risks. On one hand, it allows the model to capture qualitative changes in pipeline behavior across different eras. On the other hand, excessive flexibility risks overfitting and loss of interpretability. Our goal is therefore to introduce enough time variation to capture the dominant features of the data while maintaining parsimony and substantive meaning.

We adopt a two-part strategy. First, we restrict time variation to smooth, low-dimensional functional forms rather than allowing arbitrary year-to-year changes. Second, we use model selection criteria to assess whether the additional flexibility is warranted by the data.

To parameterize time variation, we map the calendar year $t$ to a rescaled time index $s(t)$ defined by
\begin{equation}
s(t) = \frac{t - t_{\mathrm{mid}}}{(t_{\max} - t_{\min})/2},
\end{equation}
where $t_{\min} = 1969$, $t_{\max} = 2017$, and $t_{\mathrm{mid}} = (t_{\min} + t_{\max})/2 = 1993$. 
This affine rescaling centers the observation window at zero and normalizes its length so that the endpoints lie at $s=\pm 1$, yielding a dimensionless time variable that improves numerical stability and interpretability of the time-varying parameterizations.

We model the branching fractions as smooth functions of time. Because these must lie between 0 and 1, we work on the logit scale:
\begin{equation}
\begin{aligned}
\mathrm{logit}\,\rho_{BM}(t) &= a_{BM} + b_{BM} s(t) + c_{BM} s(t)^2, \\
\mathrm{logit}\,\rho_{BP}(t) &= a_{BP} + b_{BP} s(t) + c_{BP} s(t)^2, \\
\mathrm{logit}\,\rho_{MP}(t) &= a_{MP} + b_{MP} s(t) + c_{MP} s(t)^2.
\end{aligned}
\end{equation}
As above, we do not impose a simplex constraint because these are effective routing fractions rather than cohort transition probabilities.

We allow the completion hazards $\gamma_M(t)$ and $\gamma_P(t)$ to vary over time as well. Since these are annual completion fractions, they must lie in $(0,1)$, and we again use a logistic parameterization:
\begin{equation}
\gamma_M(t) = \mathrm{logit}^{-1}(\eta_M(t)), \qquad
\gamma_P(t) = \mathrm{logit}^{-1}(\eta_P(t)),
\end{equation}
where $\eta_M(t)$ and $\eta_P(t)$ are unconstrained real-valued polynomial functions of $s(t)$.

In the most flexible specification considered, we model these as quadratic functions:
\begin{equation}
\eta_M(t) = \eta_{M0} + \eta_{M1} s(t) + \eta_{M2} s(t)^2, \qquad
\eta_P(t) = \eta_{P0} + \eta_{P1} s(t) + \eta_{P2} s(t)^2.
\end{equation}

For model comparison, we also consider an optional reduced-form forcing term added to the PhD update equation. In specifications with forcing, we include an additional exogenous input $\lambda\,p_{\mathrm{intl}}(t)$, where $p_{\mathrm{intl}}(t)$ is an external proxy series for international participation and $\lambda$ is an estimated coefficient.

Combining the time-varying branching fractions and completion hazards, the baseline model without forcing is
\begin{equation}
\begin{aligned}
M(t+1) &= M(t) + \rho_{BM}(t) \, b(t) - \gamma_M(t) \, M(t), \\
P(t+1) &= P(t) + \rho_{BP}(t) \, b(t) + \rho_{MP}(t) \, \hat{m}(t) - \gamma_P(t) \, P(t),
\end{aligned}
\end{equation}
with model-implied completion flows $\hat{m}(t) = \gamma_M(t) M(t)$ for the master's compartment and $\hat{p}(t) = \gamma_P(t) P(t)$ for the PhD compartment. 
In specifications that allow for an explicit reduced-form forcing term, the PhD update equation is augmented by an additional exogenous input,
\begin{equation}
P(t+1) = P(t) + \rho_{BP}(t) \, b(t) + \rho_{MP}(t) \, \hat{m}(t) - \gamma_P(t) \, P(t) + \lambda\, p_{\mathrm{intl}}(t),
\end{equation}
where $p_{\mathrm{intl}}(t)$ is an external proxy series for international participation (measured as the number of mathematics PhD completions awarded to temporary-visa holders, from the Survey of Earned Doctorates), and $\lambda \ge 0$ is an estimated coefficient. Because $p_{\mathrm{intl}}(t)$ is itself a completion count rather than an entry count, this term is interpreted only as a reduced-form covariate rather than as a literal demographic inflow into the PhD stock.

The fully flexible quadratic specification has 15 unknown parameters: nine governing the branching fractions and six governing the completion hazards. Simpler models with linear or constant hazards and/or branching fractions arise as constrained special cases.

From an inferential perspective, this mathematical setup defines an inverse stock--flow problem with unobserved state variables, in which identification arises from the long-run temporal structure of the observed completion flows rather than from direct observation of the underlying stocks.

\section{Estimation}\label{sec:estimation}

Given the observed bachelor's completions $b(t)$ and candidate values of the model parameters, we can iterate the discrete-time update equations forward in time to produce model-implied master's and PhD completions $\hat{m}(t)$ and $\hat{p}(t)$. The estimation problem is to find the parameter values that make these model-implied series match the observed series as closely as possible.

Because the data take the form of degree counts indexed by observation year, we work directly with the one-year updates defined above. Starting from initial conditions $M(t_{\min})$ and $P(t_{\min})$, we update the compartment populations according to
\begin{equation}
\begin{aligned}
M(t_{j+1}) &= M(t_j) + \rho_{BM}(t_j) \, b(t_j) - \gamma_M(t_j) \, M(t_j), \\
P(t_{j+1}) &= P(t_j) + \rho_{BP}(t_j) \, b(t_j) + \rho_{MP}(t_j) \, \hat{m}(t_j) - \gamma_P(t_j) \, P(t_j),
\end{aligned}
\end{equation}
and compute $\hat{m}(t_j) = \gamma_M(t_j) M(t_j)$ and $\hat{p}(t_j) = \gamma_P(t_j) P(t_j)$ at each step. These recurrences define the model exactly at the annual resolution of the data; no further numerical approximation is involved.

The initial conditions require some care. We do not observe $M(t_{\min})$ or $P(t_{\min})$ directly, so we initialize them using the observed degree counts and the fitted completion hazards:
\begin{equation}
M(t_{\min}) = \frac{m(t_{\min})}{\gamma_M(t_{\min})}, \qquad
P(t_{\min}) = \frac{p(t_{\min})}{\gamma_P(t_{\min})}.
\end{equation}
This ensures that the model-implied completions in the first observation year match the observed completions exactly, given the fitted hazard rates. The dynamics then propagate forward from this starting point. Because the system has positive completion hazards, it exhibits exponential forgetting: the influence of initial conditions decays over time, so that by the mid-period of the sample the trajectories are driven primarily by accumulated inflows rather than by the starting values. This property follows directly from the linear stock--flow structure of the model with strictly positive per-period completion hazards, which induce stable first-order dynamics. Our initialization choice may modestly affect the fitted behavior in the first few years, but our substantive conclusions focus on mid- and late-period trends, which are much less sensitive to initialization.

To compare model-implied completions to the observed degree counts, we compute residuals on the logarithmic scale. Working on the log scale has two advantages: it naturally accommodates the fact that degree counts vary over several orders of magnitude across the observation period, and it downweights years with very large counts that might otherwise dominate the fit.

For each observation year $t$, we define log-scale residuals
\begin{equation}
r_m(t) = \log m(t) - \log \hat{m}(t), \qquad r_p(t) = \log p(t) - \log \hat{p}(t).
\end{equation}
By construction, the initialization enforces $r_m(t_{\min})=r_p(t_{\min})=0$, so the first-year residuals are imposed rather than fit and do not contribute information to parameter estimation. The overall loss function is the sum of squared log-residuals across all years and both series:
\begin{equation}
\mathrm{SSE}(\theta) = \sum_t \left[ r_m(t)^2 + r_p(t)^2 \right],
\end{equation}
where $\theta$ denotes the vector of all model parameters. Because the observed degree counts in our data are strictly positive, the logarithms are well-defined. During estimation, we additionally enforce that the model-implied completions $\hat{m}(t)$ and $\hat{p}(t)$ remain strictly positive; any parameter configuration that produces non-positive fitted flows is rejected.

We estimate $\theta$ by minimizing $\mathrm{SSE}$ using a quasi-Newton optimization algorithm with the Broyden--Fletcher--Goldfarb--Shanno (BFGS) update rule. This is a standard workhorse method for smooth nonlinear least-squares problems. All free parameters enter the model through logistic or polynomial transformations, so the optimization is unconstrained in the transformed parameter space.

To quantify parameter uncertainty, we approximate the covariance matrix of the free parameters using the Hessian of the $\mathrm{SSE}$ objective evaluated at the optimum. Let $\hat{\theta}$ denote the minimizing parameter vector. We compute the observed Hessian
\begin{equation}
\widehat{H}_S = \nabla^2 \mathrm{SSE}(\theta)\big|_{\theta = \hat{\theta}}
\end{equation}
as a numerical Hessian evaluated at the optimum. We estimate the residual variance on the log scale as
\begin{equation}
\hat{\sigma}^2 = \frac{\mathrm{SSE}(\hat{\theta})}{N_{\mathrm{eff}} - k},
\end{equation}
where $k$ is the number of free parameters and $N_{\mathrm{eff}} = 96$ is the effective number of informative residual components.

Under the usual Gauss--Newton approximation for nonlinear least squares, the Hessian of $\mathrm{SSE}$ satisfies $\widehat{H}_S \approx 2 J^\top J$, where $J$ is the Jacobian of the residual vector. The familiar covariance approximation $\hat{\sigma}^2 (J^\top J)^{-1}$ thus corresponds to $2 \hat{\sigma}^2 \widehat{H}_S^{-1}$. In our implementation we therefore take
\begin{equation}
\widehat{\Sigma}_\theta = 2 \hat{\sigma}^2 \, \widehat{H}_S^{-1}
\end{equation}
as the covariance matrix for $\hat{\theta}$ and use $\mathcal{N}(\hat{\theta}, \widehat{\Sigma}_\theta)$ as a Gaussian approximation to the sampling distribution of the parameters. This curvature-based approximation is used to generate Monte Carlo confidence bands for the time-varying functions. When the numerical Hessian is ill-conditioned, we apply standard regularization and inversion safeguards to obtain a stable covariance matrix for these draws.

\section{Results}

The fully flexible quadratic specification described above has 15 free parameters. Is this level of complexity warranted by the data, or would a simpler specification suffice? To answer this question, we fit a $3 \times 3 \times 2$ grid of models that vary (i) the degree of time variation in the hazards ($\deg\gamma \in \{0,1,2\}$ for constant, linear, or quadratic in $s$), (ii) the degree of time variation in the branching fractions ($\deg\rho \in \{0,1,2\}$), and (iii) the inclusion or exclusion of an explicit international forcing term. This yields 18 candidate models in total.

For each specification, we compute the sum of squared log-residuals and evaluate both the Akaike Information Criterion (AIC) and the Bayesian Information Criterion (BIC):
\begin{equation}
\mathrm{AIC} = 2k + N \log(\mathrm{SSE}/N), \qquad \mathrm{BIC} = k \log N + N \log(\mathrm{SSE}/N),
\end{equation}
where $k$ is the number of parameters and $N$ is the total number of observations (49 master's observations plus 49 PhD observations, so $N = 98$). Lower values indicate better fit after penalizing for complexity. Because our initialization forces the first-year residuals to be exactly zero for both series, one could alternatively view the effective sample size as $N_{\mathrm{eff}} = 96$ rather than $N = 98$. Using $N_{\mathrm{eff}}$ would slightly change the absolute AIC and BIC values (and the BIC penalty via $\log N$), but in our setting it does not affect the qualitative ranking of models or the substantive conclusions. Table~\ref{tab:model_comparison} reports the results, with quantities rounded for readability.

\begin{table}[ht]
\caption{Model comparison across 18 specifications. Columns indicate the polynomial degree of the hazards ($\deg\gamma$), the degree of the branching fractions ($\deg\rho$), whether an explicit international forcing term is included, the number of free parameters $k$, the sum of squared log-residuals (SSE), and information criteria. $\Delta$AIC and $\Delta$BIC are relative to the best model (Model 1: $\deg\gamma=2$, $\deg\rho=2$, no forcing).}
\label{tab:model_comparison}
\centering
\begin{tabular}{cccccrrrrr}
\toprule
Model & $\deg\gamma$ & $\deg\rho$ & Inflow & $k$ & SSE & AIC & $\Delta$AIC & BIC & $\Delta$BIC \\
\midrule
1  & 2 & 2 & N & 15 & 0.129 & $-$619.7 & 0.0  & $-$581.0 & 0.0 \\
2  & 2 & 2 & Y & 16 & 0.129 & $-$617.7 & 2.0  & $-$576.4 & 4.6 \\
3  & 1 & 2 & N & 13 & 0.170 & $-$597.2 & 22.5 & $-$563.6 & 17.3 \\
4  & 1 & 2 & Y & 14 & 0.170 & $-$594.8 & 24.9 & $-$558.6 & 22.4 \\
5  & 2 & 1 & N & 12 & 0.229 & $-$569.8 & 50.0 & $-$538.7 & 42.2 \\
6  & 2 & 1 & Y & 13 & 0.229 & $-$567.8 & 52.0 & $-$534.1 & 46.8 \\
7  & 1 & 1 & N & 10 & 0.297 & $-$548.4 & 71.3 & $-$522.6 & 58.4 \\
8  & 1 & 1 & Y & 11 & 0.297 & $-$546.4 & 73.3 & $-$518.0 & 62.9 \\
9  & 0 & 2 & N & 11 & 0.299 & $-$545.7 & 74.1 & $-$517.2 & 63.7 \\
10 & 0 & 2 & Y & 12 & 0.299 & $-$543.7 & 76.1 & $-$512.7 & 68.3 \\
11 & 0 & 1 & N &  8 & 0.656 & $-$474.6 & 145.1 & $-$454.0 & 127.0 \\
12 & 0 & 1 & Y &  9 & 0.656 & $-$472.6 & 147.1 & $-$449.4 & 131.6 \\
13 & 2 & 0 & N &  9 & 0.684 & $-$468.5 & 151.2 & $-$445.3 & 135.7 \\
14 & 2 & 0 & Y & 10 & 0.684 & $-$466.5 & 153.2 & $-$440.7 & 140.3 \\
15 & 0 & 0 & N &  5 & 6.028 & $-$263.3 & 356.5 & $-$250.3 & 330.6 \\
16 & 0 & 0 & Y &  6 & 6.028 & $-$261.3 & 358.5 & $-$245.8 & 335.2 \\
17 & 1 & 0 & N &  7 & 6.029 & $-$259.3 & 360.5 & $-$241.2 & 339.8 \\
18 & 1 & 0 & Y &  8 & 6.029 & $-$257.3 & 362.5 & $-$236.6 & 344.4 \\
\bottomrule
\end{tabular}
\end{table}

The model grid yields three concrete conclusions. First, both the branching fractions and the completion hazards require quadratic time variation: any reduction to $\deg\rho<2$ or $\deg\gamma<2$ worsens fit enough that the penalized criteria degrade substantially (the best model with $\deg\gamma<2$ is worse by $\Delta\mathrm{AIC}=22.5$ and $\Delta\mathrm{BIC}=17.3$). Second, adding an explicit international forcing term provides no improvement for the quadratic--quadratic specification and worsens both AIC (by about 2 points) and BIC (by about 5 points). Third, the overall best model by both AIC and BIC is the 15-parameter specification with quadratic branching fractions, quadratic completion hazards, and no explicit forcing (Model 1 in Table~\ref{tab:model_comparison}).

\begin{figure}[ht]
\centering
\includegraphics[width=\textwidth]{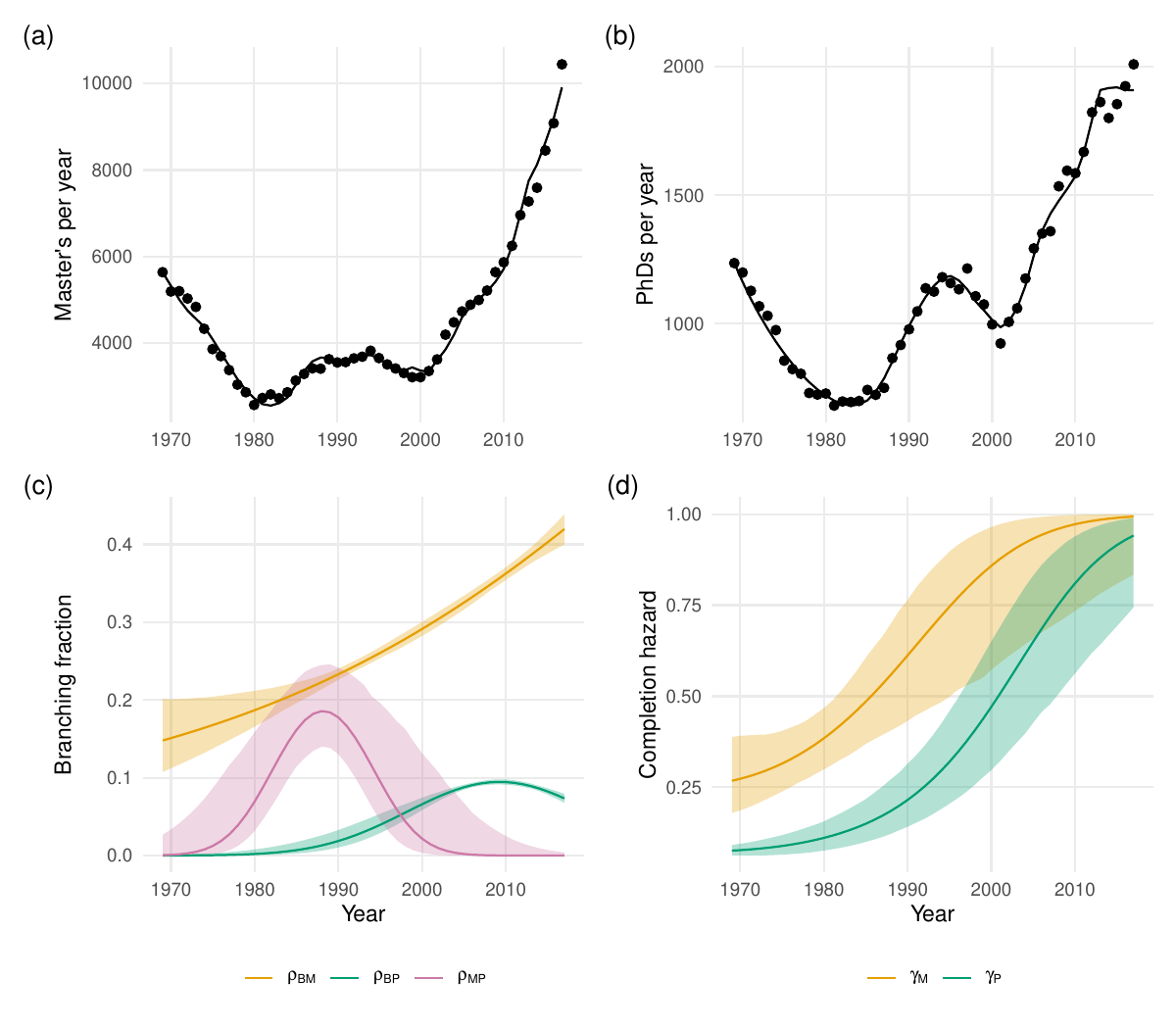}
\caption{Fitted results from the preferred model (quadratic branching fractions with quadratic time-varying completion hazards, no explicit forcing). Panels (a) and (b) show observed degree counts (points) and model-implied fits (lines) for master's and PhD completions, respectively, over the 1969--2017 window. Panel (c) displays the estimated branching fractions over time, with shaded 95\% Monte Carlo confidence bands constructed from parameter draws based on the Hessian: $\rho_{BM}(t)$ (bachelor's to master's), $\rho_{BP}(t)$ (bachelor's directly to PhD), and $\rho_{MP}(t)$ (master's to PhD). Panel (d) shows the estimated completion hazards $\gamma_M(t)$ and $\gamma_P(t)$ with analogous 95\% bands.}
\label{fig:results}
\end{figure}

We therefore adopt this quadratic--quadratic, no-forcing specification as our preferred model. Figure~\ref{fig:results} displays the fitted results along with 95\% confidence bands computed via Monte Carlo simulation from the Hessian-based parameter covariance matrix.

To interpret the quality of fit, we convert the sum of squared log-residuals into a more intuitive measure. The preferred model achieves a sum of squared log-residuals of $\mathrm{SSE} \approx 0.129$ across $N = 98$ observations. The root-mean-squared error (RMSE) on the log scale is therefore
\begin{equation}
\sqrt{0.129/98} \approx 0.036.
\end{equation}
For small residuals, a log-residual $r = \log(\text{observed}/\text{fitted})$ is approximately equal to the relative error $(\text{observed} - \text{fitted})/\text{fitted}$, since $\log(1 + x) \approx x$ when $x$ is small. Thus an RMSE of about 0.036 on the log scale corresponds to typical multiplicative errors of about $e^{0.036}-1 \approx 3.7\%$ in the fitted degree counts. The model captures the major qualitative features of both the master's and PhD degree series, including periods of decline, recovery, and rapid growth.

The confidence bands in Figure~\ref{fig:results} quantify uncertainty arising from the fact that multiple parameter configurations can fit the observed degree series almost equally well. We approximate the sampling distribution of the free parameters using the Gaussian approximation described above, with covariance $\widehat{\Sigma}_\theta = 2 \hat{\sigma}^2 \widehat{H}_S^{-1}$, and then draw many parameter vectors from the resulting multivariate normal distribution. For each draw we recompute $\rho_{BM}(t)$, $\rho_{BP}(t)$, $\rho_{MP}(t)$, $\gamma_M(t)$, and $\gamma_P(t)$; the shaded regions show the 2.5th and 97.5th percentiles across draws at each observation year.

Because the observed data constrain only completion flows rather than the latent stocks themselves, branching fractions---which govern how observed inflows are partitioned across pathways---are more directly identified than completion hazards, which trade off against unobserved stock sizes to produce the same observed exits. The branching fraction bands are relatively narrow, indicating that the data strongly constrain how the bachelor's and master's completion flows feed forward into downstream compartments. By contrast, the hazard bands in panel (d) of Figure~\ref{fig:results} are noticeably wider. This is not a numerical artifact but a structural feature of the model: we observe only the completion flows $\gamma_M(t)M(t)$ and $\gamma_P(t)P(t)$, not the latent compartment populations themselves. As a result, many different combinations of hazards and latent stocks can produce very similar completion trajectories, so the hazards are inherently less well identified than the branching fractions.

These narrow bands are \emph{marginal} confidence bands for each branching fraction considered individually. However, joint identification of the PhD inflow decomposition is weaker, as we now discuss. Relatedly, because the data constrain only the aggregate PhD completion flow $\hat p(t)=\gamma_P(t)P(t)$ (and not PhD entry counts or the latent stock $P(t)$), the model's late-period decomposition of PhD inflows into the terms $\rho_{BP}(t)b(t)$ and $\rho_{MP}(t)\hat m(t)$ is not separately well identified from $\gamma_P(t)$ and the latent stock $P(t)$, and should be read as a reduced-form accounting device rather than as a literal statement about where PhD entrants come from. Two related but distinct identification issues are at work here. First, because only the product $\hat p(t)=\gamma_P(t)P(t)$ is constrained by the data, the PhD hazard $\gamma_P(t)$, the latent stock $P(t)$, and the total inflow into $P$ are jointly non-separable. Second, even conditional on a given total inflow, uncertainty in $\hat m(t)$ and joint covariance between $\rho_{BP}(t)$ and $\rho_{MP}(t)$ can leave the inflow partition ambiguous even when the marginal confidence bands for the individual branching fractions are relatively narrow. The point estimates in panel (d) therefore represent the most likely hazard trajectories consistent with the data and model, but the confidence bands remind us that a range of trajectories is observationally similar. By the end of the observation window, for example, the 95\% band for $\gamma_M(t)$ spans a nontrivial range, so the magnitude of the apparent speeding-up of master's completion should be interpreted as suggestive rather than precisely estimated. Moreover, because hazards are constrained to lie in $(0,1)$, late-period point estimates that approach the upper boundary should be interpreted as indicating near-maximal effective turnover given the latent stock scaling implied by the model, not as evidence that typical individual time-to-degree approaches one year.

\begin{figure}[ht]
\centering
\includegraphics[width=\textwidth]{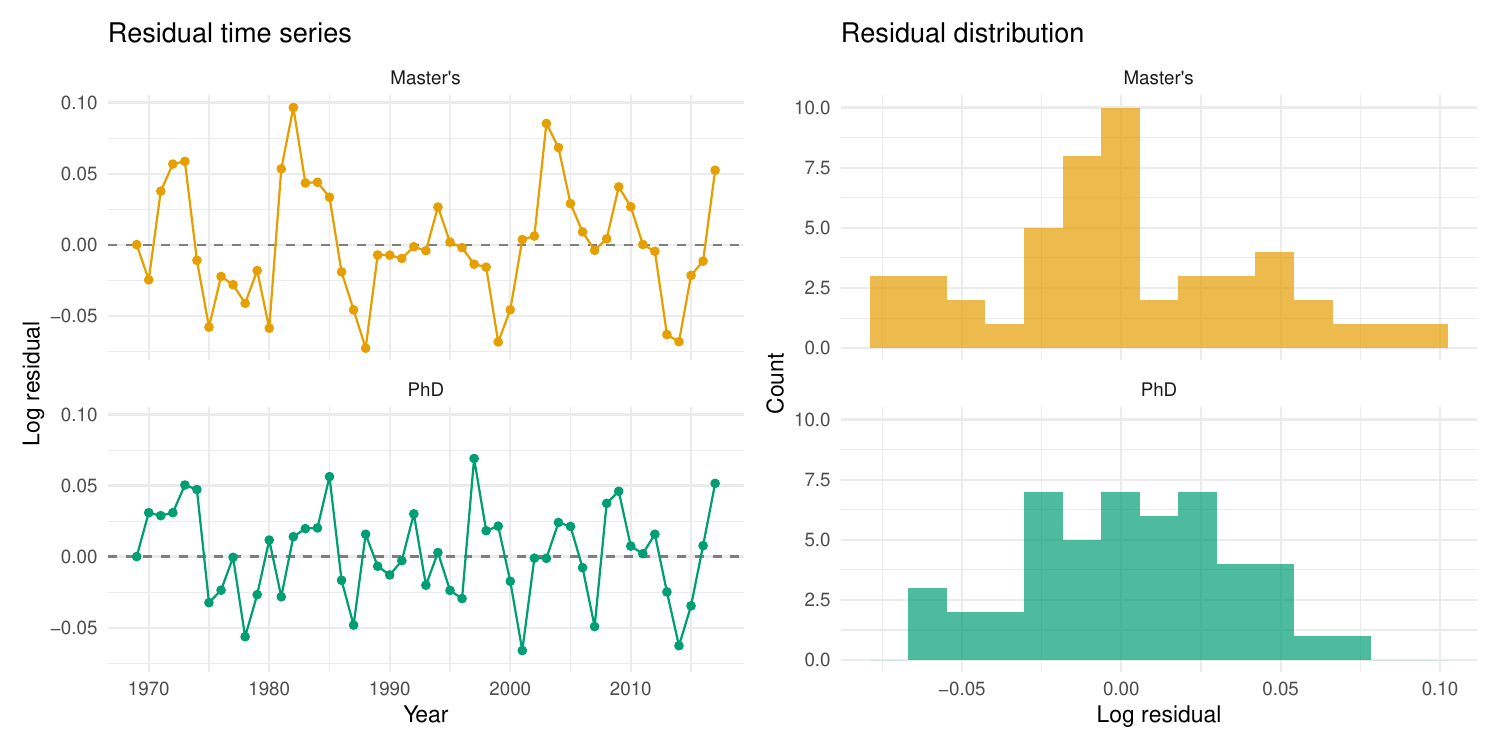}
\caption{Residual diagnostics for the preferred model. Left: log-residuals for master's and PhD completions over time. Right: histograms of log-residuals. Residuals are centered near zero and mostly lie within approximately $\pm 0.1$, consistent with typical relative errors of a few percent.}
\label{fig:resid}
\end{figure}

We present residual diagnostics in Figure~\ref{fig:resid}. The left panel shows log-residuals for master's and PhD completions over time; the right panel shows their empirical distributions. Residuals are centered near zero. Master's residuals are modestly positive in several late-period years, consistent with the possibility that even the quadratic specification slightly underestimates the very recent surge in master's completions, but we do not see a large, persistent trend. The residual distributions are approximately symmetric and concentrated within $\pm 0.1$ on the log scale, consistent with the RMSE calculation above.

Finally, we conduct two targeted robustness checks aimed at assessing whether the preferred 15-parameter structure is unduly dependent on the earliest observations or on short-run idiosyncrasies.

The first check evaluates sensitivity to the early period and initialization. If the estimated trajectories were heavily driven by early-period leverage or by our initialization convention, we would expect fit quality to deteriorate materially when the earliest years are excluded. To test this, we refit the preferred specification after discarding the first 5, 10, and 15 years of the sample (i.e., re-estimating on 1974--2017, 1979--2017, and 1984--2017). All three truncated-sample refits converged and achieved pooled log-RMSE (computed over the truncated window) in the narrow range 0.031--0.032, comparable to or slightly lower than the full-sample value of 0.036 (Table~\ref{tab:robust_startyear}). On the multiplicative scale, this corresponds to typical errors of about $e^{0.032}-1 \approx 3.3\%$. The stability of fit quality across these truncations indicates that the preferred fit is not driven by the earliest years.

\begin{table}[ht]
\caption{Start-year sensitivity for the preferred specification using systematic 5-, 10-, and 15-year truncations of the early sample. For each start year, we refit the preferred model on the truncated window (start year through 2017). Pooled log-RMSE is computed as $\sqrt{\mathrm{SSE}/N}$ with $N = 2 \times (\text{number of years in the truncated window})$, reflecting pooled evaluation across the master's and PhD series.}
\label{tab:robust_startyear}
\centering
\begin{tabular*}{\textwidth}{@{\extracolsep{\fill}}ccccc@{}}
\toprule
Start year & Converged & $k$ & SSE & Pooled log-RMSE \\
\midrule
1974 & True & 15 & 0.0899 & 0.0320 \\
1979 & True & 15 & 0.0804 & 0.0321 \\
1984 & True & 15 & 0.0664 & 0.0312 \\
\bottomrule
\end{tabular*}
\end{table}

The second check assesses short-horizon stability via rolling-origin one-step-ahead hindcasts. If the fitted model were merely absorbing year-to-year noise rather than capturing stable structure, we would expect poor performance when asked to predict the next year's completions from parameters estimated on earlier data. To test this, for each cutoff year $T \in \{1990, 1995, 2000, 2005, 2010, 2015\}$ we refit the preferred model using only data through year $T$ and then predict degree completions in year $T+1$. We compute the one-step log-RMSE separately for each series:
\begin{equation}
\mathrm{RMSE}_{\mathrm{1step}} = \sqrt{\frac{1}{J} \sum_{j=1}^{J} \bigl(\log y_{T_j+1} - \log \hat y_{T_j+1}^{(T_j)}\bigr)^2},
\end{equation}
where $J=6$ is the number of cutoffs and $\hat y_{T_j+1}^{(T_j)}$ denotes the model prediction for year $T_j+1$ based on parameters estimated using data through year $T_j$. This rolling-origin evaluation yields $\mathrm{RMSE}_{\mathrm{1step}} = 0.0723$ for master's completions and $0.0735$ for PhD completions (pooled across both series: $0.0729$). On the multiplicative scale, these correspond to typical one-step errors of about $e^{0.073}-1 \approx 7.6\%$. This is larger than the full-sample in-sample RMSE of approximately 0.036 (about $e^{0.036}-1 \approx 3.7\%$), as expected for an out-of-sample check, but it remains modest relative to the large secular changes in degree production over the observation window.

\section{Discussion}

The fitted model reveals several noteworthy patterns in the evolution of the mathematics training pipeline.

First, the way students enter PhD programs has rebalanced over time. 
Panel (c) of Figure~\ref{fig:results} shows that $\rho_{BM}(t)$---the bachelor's-to-master's branching fraction---rises steadily from roughly 0.2 to above 0.35 over the sample period, indicating that a growing share of bachelor's recipients continue into mathematics master's programs.
The branching fraction $\rho_{BP}(t)$---direct bachelor's-to-PhD entry---remains small throughout (well below 0.1), with only modest changes over time. 
By contrast, $\rho_{MP}(t)$---the fraction of master's completions that contribute to PhD entry---is strongly time-varying: it increases from low levels in the early 1970s to a peak in the late 1980s and early 1990s before declining again in the 2000s and 2010s. 
Together, these patterns in the inferred branching fractions suggest that the master's degree became increasingly central to the PhD pipeline through the mid-period of the sample, but that its role as the dominant gateway into doctoral study has softened in the sense that the conditional master's-to-PhD fraction $\rho_{MP}(t)$ declined, even as the absolute contribution of the master's pathway may have continued to grow as master's completions increased.

Second, panel (d) shows that the point estimates for both completion hazards rise over time, particularly for master's programs. 
For master's programs, $\gamma_M(t)$ increases substantially, implying shorter effective durations in the master's compartment under the model's reduced-form turnover mechanism. For PhD programs, $\gamma_P(t)$ also rises, though more modestly. 
Interpreted as effective turnover, these trajectories suggest that the pipeline has become faster-moving over time, with shorter effective lags between entering a compartment and completing the degree. 
At the same time, the wide confidence bands for the hazards underscore that these trends are estimated with more uncertainty than the branching fractions: alternative hazard trajectories within the bands are also consistent with the observed degree series, and the range of plausible values for late-period $\gamma_M(t)$ spans roughly a factor of four.

These mechanical changes in routing behavior and effective turnover are consistent with several qualitative developments in U.S. graduate education over this period.
The growth of professional and terminal master's programs, shifts in funding and assistantship patterns, changes in faculty hiring and job-market signals, and evolving visa and immigration policies may all have contributed to the strengthening of the bachelor's-to-master's pathway and to changes in the master's-to-PhD channel.
Our aggregate model is not designed to test these mechanisms individually, but the fitted trajectories provide quantitative targets that future work with richer microdata could attempt to explain.

The model treats the observed bachelor's, master's, and PhD degree counts as aggregate flows and does not distinguish between domestic and international students. Over the past several decades, however, the composition of PhD recipients has shifted markedly: data from the Survey of Earned Doctorates (SED) indicate a pronounced increase in the number of doctorates awarded to temporary-visa holders in mathematics and statistics \cite{NCSES2024DoctorateRecipients2023}.

One might therefore wonder whether incorporating an explicit reduced-form proxy for international participation would improve the model.
To investigate this, we considered specifications that add an exogenous forcing term to the PhD update equation proportional to the SED-based count of temporary-visa PhD completions, $\lambda\,p_{\mathrm{intl}}(t)$. 
Because $p_{\mathrm{intl}}(t)$ is a completion count rather than an entry count, this term is interpreted only as a reduced-form proxy covariate rather than as a literal demographic inflow into the PhD stock. 
Table~\ref{tab:model_comparison} shows that including this term does not improve the model: both AIC and BIC worsen for the quadratic--quadratic specification, and the estimated $\lambda$ is driven close to zero across specifications.
Within this particular reduced-form test, the time-varying branching fractions already capture the long-run structural changes that the temporary-visa series would otherwise proxy.

This finding reinforces the reduced-form interpretation of our model. 
The fitted branching fractions and completion hazards are effective parameters that summarize how bachelor's and master's completion flows contribute to downstream degree production in the aggregate. 
They incorporate the net effects of both domestic and international students, of various entry pathways, and of other sources of heterogeneity that we do not model explicitly.
The model provides a coherent dynamical framework for summarizing long-run evolution in the mathematics training pipeline but it is not designed to recover separate domestic and international sub-pipelines or to measure individual-level transition probabilities. 
Because the model is observational and reduced-form, the branching fractions should not be interpreted as causal responses to changes in $b(t)$ or $m(t)$.
The current specification also treats inflows as contemporaneous; in reality, there are lags between degree completion at one stage and entry at the next, and multi-year enrollment durations within each stage.
Extending the framework to incorporate explicit lags or additional compartments would be a natural next step, but would require data beyond the aggregate flows analyzed here.

More broadly, the analysis illustrates how linked degree-completion time series can be recast as a dynamically coherent system even when intermediate stocks are entirely unobserved.
This perspective is potentially useful in other settings where only exits are observed but the underlying stock--flow dynamics are of substantive interest.

\section{Conclusion}

We have developed a minimal compartmental model that links the annual production of bachelor's, master's, and PhD degrees in mathematics over nearly five decades.
The model introduces latent compartments representing unobserved graduate populations and relates these populations to observed degree flows through time-varying branching fractions and completion hazards. 
The model is purely descriptive: it summarizes aggregate patterns in the data within a reduced-form framework without claiming causal identification.

The key features of the model can be summarized as follows. 
First, the model takes bachelor's degree completions as an exogenous input and produces master's and PhD completions as outputs, capturing the idea that graduate degree production depends on the accumulated stock of students flowing through the pipeline. 
Second, the model parameters are allowed to vary smoothly over time, which turns out to be essential for fitting the data: both the branching fractions and the completion hazards exhibit substantial trends across the observation period.
Third, model selection criteria confirm that a specification with quadratic time variation in both branching fractions and completion hazards, and no explicit international forcing term, is strongly preferred.
Both AIC and BIC favor this no-forcing specification over the otherwise identical quadratic--quadratic model that includes the forcing term, indicating that the forcing term fails to improve fit sufficiently to justify its additional parameter.

What does the model reveal that simple inspection of the raw degree series would not? 
Three findings stand out.
First, the structure of PhD entry has rebalanced over time: the bachelor's-to-master's pathway strengthens, the master's-to-PhD channel intensifies through the 1980s before receding, and direct bachelor's-to-PhD entry remains a relatively small but persistent route, as captured by the model's inferred branching fractions.
Second, the model's estimated completion hazards increase over time, suggesting that the pipeline has become faster-moving in an effective, reduced-form sense, although these hazard trends are estimated with considerably more uncertainty than the branching behavior. 
Third, these patterns emerge as smooth, gradual evolution rather than abrupt regime changes, with the reorganization of the pipeline unfolding continuously over decades.

Several limitations should be kept in mind when interpreting the results. 
The fitted parameters are effective aggregate quantities, not direct measurements of individual-level transition probabilities or completion times. 
The model does not distinguish between domestic and international students, between different types of master's degrees, or between different institutional pathways through the pipeline, and it does not explicitly model lags between degree completion at one stage and entry at the next.
These simplifications are necessary given the aggregate nature of the available data, and they mean that the model is best understood as a reduced-form summary of pipeline dynamics rather than a structural model of individual decisions.

Despite these limitations, the two-compartment framework provides a useful lens for understanding long-run evolution in mathematics degree production.
It quantifies how changes in bachelor's output are reflected in downstream master's and PhD degree production over time, how effective turnover has changed, and how the relative contribution of different entry pathways has shifted across decades. 
More broadly, it provides a general approach for studying stock--flow systems observed only through exits, offering a way to reconstruct and summarize long-run routing and turnover dynamics from linked flow time series when intermediate stocks are unobserved. 
Future work could extend this framework to incorporate additional compartments, explicit lag structures, or disaggregated domestic and international subpopulations, using richer data sources to probe the mechanisms that underlie the aggregate patterns documented here.

%%======================================================================
%% BACK MATTER
%%======================================================================
\backmatter

\ifpreprint
  %%--------------------------------------------------------------------
  %% Statements and Declarations (included in preprint)
  %%--------------------------------------------------------------------
  \bmhead{Acknowledgments}
  This work was conceived at the American Institute of Mathematics (AIM) workshop ``MetaMath: Modeling the mathematical sciences community using mathematics, statistics, and data science,'' held December 8--12, 2025 in Pasadena, California. The authors thank AIM for its hospitality and the workshop participants for valuable discussions. AIM is supported by the National Science Foundation.

  \bmhead{Funding}
  No funding was received for conducting this study.

  \bmhead{Competing Interests}
  The authors declare no competing interests relevant to the content of this article.

  \bmhead{Data Availability}
  The degree completion data analyzed in this study are publicly available from the National Center for Education Statistics (NCES), U.S.\ Department of Education, at \url{https://nces.ed.gov/}.

  \bmhead{Code Availability}
  The analysis code is available at \url{https://doi.org/10.5281/zenodo.18222401}.

  \bmhead{Disclaimer}  The views expressed in this article are those of the author(s) and do not necessarily reflect the official policy or position of the United States Air Force Academy, the Air Force, the Department of Defense, or the U.S. Government. [PA number forthcoming]    

  \bmhead{Author Contributions}
  C.M.T.\ conceived the study, developed the latent-stock compartmental framework, wrote the software, performed the analysis, and wrote the manuscript. O.B., R.B., D.G., M.H., O.I., D.N., and A.R.V.-M.\ contributed to conceptual discussions at the AIM workshop and provided feedback on the manuscript. All authors reviewed and approved the final manuscript.
  C.M.T. is the lead author; all other authors are listed alphabetically by last name.
\else
  %%--------------------------------------------------------------------
  %% For double-blind review: declarations go on separate title page
  %%--------------------------------------------------------------------
\fi

%%----------------------------------------------------------------------
%% References
%%----------------------------------------------------------------------
\bibliography{bibliography}

\end{document}